\documentclass{article}

\usepackage{graphicx,etoolbox}

\usepackage{hyperref}

\usepackage{hyperref}

\usepackage{amssymb}

\newcommand{\newc}{\newcommand}

\newc{\eqnoset}{\setcounter{equation}{0}}

\newcommand{\mref}[1]{(\ref{#1})}

\newcommand{\refrem}[1]{Remark~\ref{#1}}
\newcommand{\reftheo}[1]{Theorem~\ref{#1}}

\newcommand{\refcoro}[1]{Corollary~\ref{#1}}

\newcommand{\refsec}[1]{Section~\ref{#1}}

\newcommand{\beq}{\begin{equation}}
	\newcommand{\eeq}{\end{equation}}
\newcommand{\beqno}[1]{\begin{equation}\label{#1}}
	
	\newcommand{\barr}{\begin{array}}
		\newcommand{\earr}{\end{array}}
	
	\newc{\bearr}{\begin{eqnarray*}}
		\newc{\eearr}{\end{eqnarray*}}
	
	\newc{\bearrno}[1]{\begin{eqnarray}\label{#1}}
		\newc{\eearrno}{\end{eqnarray}}
	
	\newc{\non}{\nonumber}
	\newc{\nol}{\nonumber\nl}
	
	\newcommand{\bdes}{\begin{description}}
		\newcommand{\edes}{\end{description}}
	\newc{\benu}{\begin{enumerate}}
		\newc{\eenu}{\end{enumerate}}
	\newc{\btab}{\begin{tabular}}
		\newc{\etab}{\end{tabular}}

	\newtheorem{theorem}{Theorem}[section]
	\newtheorem{defi}[theorem]{Definition}
	\newtheorem{lemma}[theorem]{Lemma}
	\newtheorem{rem}[theorem]{Remark}
	\newtheorem{exam}[theorem]{Example}
	\newtheorem{propo}[theorem]{Proposition}
	\newtheorem{corol}[theorem]{Corollary}
	\newtheorem{conj}[theorem]{Conjecture}

	\newcommand{\btheo}[1]{\begin{theorem}\label{#1}}
		\newc{\brem}[1]{\begin{rem}\label{#1}\em}
			\newc{\bexam}[1]{\begin{exam}\label{#1}\em}
				\newc{\bdefi}[1]{\begin{defi}\label{#1}}
					\newcommand{\blemm}[1]{\begin{lemma}\label{#1}}
						\newcommand{\bprop}[1]{\begin{propo}\label{#1}}
							\newcommand{\bcoro}[1]{\begin{corol}\label{#1}}
								\newcommand{\btheoc}[1]{\begin{conj}\label{#1}}
									\newcommand{\etheo}{\end{theorem}}
								\newc{\etheoc}{\end{conj}}
							\newcommand{\elemm}{\end{lemma}}
						\newcommand{\eprop}{\end{propo}}
					\newcommand{\ecoro}{\end{corol}}
				\newc{\erem}{\end{rem}}
			\newc{\eexam}{\end{exam}}
		\newc{\edefi}{\end{defi}}
	
	\newc{\rmk}[1]{{\bf REMARK #1: }}
	\newc{\DN}[1]{{\bf DEFINITION #1: }}
	
	\newcommand{\bproof}{{\bf Proof:~~}}
	\newc{\eproof}{{\vrule height8pt width5pt depth0pt}\vspace{3mm}}
	
	\newc{\bfrac}[2]{\dspl{\frac{#1}{#2}}}

	\newc{\nid}{\noindent}
	
	\newcommand{\dspl}{\displaystyle}
	\newc{\grad}{\nabla}
	\newc{\Div}{\mbox{div}}
	\newc{\pdt}[1]{\dspl{\frac{\partial{#1}}{\partial t}}}
	\newc{\pdn}[1]{\dspl{\frac{\partial{#1}}{\partial \nu}}}
	\newc{\pdNi}[1]{\dspl{\frac{\partial{#1}}{\partial \mathcal{N}_i}}}
	\newc{\pD}[2]{\dspl{\frac{\partial{#1}}{\partial #2}}}
	\newc{\dt}{\dspl{\frac{d}{dt}}}
	\newc{\bdry}[1]{\mbox{$\partial #1$}}
	\newc{\sgn}{\mbox{sign}}
	
	\newc{\Hess}[1]{\frac{\partial^2 #1}{\pdh z_i \pdh z_j}}
	\newc{\hess}[1]{\partial^2 #1/\pdh z_i \pdh z_j}

	\newc{\ag}{\alpha}
	\newc{\bg}{\beta}
	\newc{\cg}{\gamma}\newc{\Cg}{\Gamma}
	\newc{\dg}{\delta}\newc{\Dg}{\Delta}
	\newc{\eg}{\varepsilon}
	\newc{\zg}{\zeta}
	\newc{\thg}{\theta}
	\newc{\llg}{\lambda}\newc{\LLg}{\Lambda}
	\newc{\kg}{\kappa}
	\newc{\rg}{\rho}
	\newc{\sg}{\sigma}\newc{\Sg}{\Sigma}
	\newc{\tg}{\tau}
	\newc{\fg}{\phi}\newc{\Fg}{\Phi}
	\newc{\vfg}{\varphi}
	\newc{\og}{\omega}\newc{\Og}{\Omega}
	\newc{\pdh}{\partial}
	
	\newc{\ccG}{{\cal G}}

	\newc{\ii}[1]{\int_{#1}}
	\newc{\iidx}[2]{{\dspl\int_{#1}~#2~dx}}
	\newc{\bii}[1]{{\dspl \ii{#1} }}
	\newc{\biii}[2]{{\dspl \iii{#1}{#2} }}
	\newc{\su}[2]{\sum_{#1}^{#2}}
	\newc{\bsu}[2]{{\dspl \su{#1}{#2} }}

	\newc{\biiom}[1]{{\dspl\int_{\bdrom}~ #1 ~d\sg}}
	\newc{\io}[1]{{\dspl\int_{\Og}~ #1 ~dx}}
	\newc{\bio}[1]{{\dspl\int_{\bdrom}~ #1 ~d\sg}}
	\newc{\bsir}{\bsu{i=1}{r}}
	\newc{\bsim}{\bsu{i=1}{m}}
	
	\newc{\iibr}[2]{\iidx{\bprw{#1}}{#2}}
	\newc{\Intbr}[1]{\iibr{R}{#1}}
	\newc{\intbr}[1]{\iibr{\rg}{#1}}
	\newc{\intt}[3]{\int_{#1}^{#2}\int_\Og~#3~dxdt}
	
	\newc{\itQ}[2]{\dspl{\int\hspace{-2.5mm}\int_{#1}~#2~dz}}
	\newc{\mitQ}[2]{\dspl{\rule[1mm]{4mm}{.3mm}\hspace{-5.3mm}\int\hspace{-2.5mm}\int_{#1}~#2~dz}}
	\newc{\mitQQ}[3]{\dspl{\rule[1mm]{4mm}{.3mm}\hspace{-5.3mm}\int\hspace{-2.5mm}\int_{#1}~#2~#3}}
	
	\newc{\mitx}[2]{\dspl{\rule[1mm]{3mm}{.3mm}\hspace{-4mm}\int_{#1}~#2~dx}}
	\newc{\mitmu}[2]{\dspl{\rule[1mm]{3mm}{.3mm}\hspace{-4mm}\int_{#1}~#2~d\mu}}
	\newc{\iidmu}[2]{\iidx{#1}{#2}}
	
	\newc{\iidm}[3]{{\dspl\int_{#1}~#2~d #3}}
	
	\newc{\itQmu}[2]{\dspl{\int\hspace{-2.5mm}\int_{#1}~#2~d\mu}}
	\newc{\mitQmu}[2]{\dspl{\rule[1mm]{4mm}{.3mm}\hspace{-5.3mm}\int\hspace{-2.5mm}\int_{#1}~#2~d\mu}}

	\newc{\mitQq}[2]{\dspl{\rule[1mm]{4mm}{.3mm}\hspace{-5.3mm}\int\hspace{-2.5mm}\int_{#1}~#2~d\bar{z}}}
	\newc{\itQq}[2]{\dspl{\int\hspace{-2.5mm}\int_{#1}~#2~d\bar{z}}}

	\newc{\pder}[2]{\dspl{\frac{\partial #1}{\partial #2}}}

	\newc{\bdrom}{\bdry{\Og}}

	\newc{\bilhom}{\mbox{Bil}(\mbox{Hom}(\RR^{nm},\RR^{nm}))}
	\newc{\VV}[1]{{V(Q_{#1})}}
	
	\newc{\ccA}{{\mathcal A}}
	\newc{\ccB}{{\mathcal B}}
	\newc{\ccC}{{\mathcal C}}
	\newc{\ccD}{{\mathcal D}}
	\newc{\ccE}{{\mathcal E}}
	\newc{\ccH}{\mathcal{H}}
	\newc{\ccF}{\mathcal{F}}
	\newc{\ccI}{{\mathcal I}}
	\newc{\ccJ}{{\mathcal J}}
	\newc{\ccK}{{\mathcal K}}
	\newc{\ccP}{{\mathcal P}}
	\newc{\ccQ}{{\mathcal Q}}
	\newc{\ccR}{{\mathcal R}}
	\newc{\ccS}{{\mathcal S}}
	\newc{\ccT}{{\mathcal T}}
	\newc{\ccX}{{\mathcal X}}
	\newc{\ccY}{{\mathcal Y}}
	\newc{\ccZ}{{\mathcal Z}}
	
	\newc{\bb}[1]{{\mathbf #1}}
	
	\newc{\myprod}[1]{\langle #1 \rangle}
	\newc{\mypar}[1]{\left( #1 \right)}
	
	\newc{\BLLg}{\mathbf{\LLg}}
	
	\newc{\mA}{\mathbf{A}}
	\newc{\mB}{\mathbf{B}}
	\newc{\mC}{\mathbf{C}}
	\newc{\mD}{\mathbf{D}}
	\newc{\mE}{\mathbf{E}}
	\newc{\mF}{\mathbf{F}}
	\newc{\mJ}{\mathbf{J}}
	\newc{\mG}{\mathbf{G}}
	\newc{\mP}{\mathbf{P}}
	\newc{\mR}{\mathbf{R}}
	\newc{\mQ}{\mathbf{Q}}
	\newc{\mX}{\mathbf{X}}
	\newc{\muu}{\mathbf{u}}
	\newc{\mvv}{\mathbf{v}}
	
	\newc{\mllg}{\mathbb{\lambda}}
	\newc{\mLLg}{\mathbf{\LLg}}

	\newc{\lspn}[2]{\mbox{$\| #1\|_{\Lsp{#2}}$}}
	\newc{\Lpn}[2]{\mbox{$\| #1\|_{#2}$}}
	\newc{\Hn}[1]{\mbox{$\| #1\|_{H^1(\Og)}$}}

	\newc{\mynorm}[2]{\| #1\|_{#2}}

	\newcommand{\RR}{{\rm I\kern -1.6pt{\rm R}}}

	\newc{\itQQ}[2]{\dspl{\int_{#1}#2\,dz}}
	\newc{\mmitQQ}[2]{\dspl{\rule[1mm]{4mm}{.3mm}\hspace{-4.3mm}\int_{#1}~#2~dz}}
	\newc{\MmitQQ}[2]{\dspl{\rule[1mm]{4mm}{.3mm}\hspace{-4.3mm}\int_{#1}~#2~d\mu}}

	\newc{\MUmitQQ}[3]{\dspl{\rule[1mm]{4mm}{.3mm}\hspace{-4.3mm}\int_{#1}~#2~d#3}}
	\newc{\MUitQQ}[3]{\dspl{\int_{#1}~#2~d#3}}

	\newc{\mccP}{\mathbb{P}}
	\newc{\mccK}{\mathbb{K}}
	
	\newc{\DKTmU}{\mccK(U)}
	\newc{\DKTmUold}{(K_U(U)^{-1})^T}
	
	\newc{\myPi}{\mathbf{W}}
	\newc{\myIbar}{\bar{\ccI}_1}
	\newc{\myIhat}{\hat{\ccI}_1}
	\newc{\myIbreve}{\breve{\ccI}_0}
	
	\newc{\mmk}{\mathbf{k}}

	\newcommand{\ma}{\mathbf{a}}

	\newc{\mfu}{\mathbf{f_u}}
	\newc{\mh}{\mathbf{h}}
	\newc{\mb}{\mathbf{b}}
	\newc{\mf}{\mathbf{f}}

	\newc{\twomatrix}[1]{\left[\barr{cc}#1\earr\right]}
	\newc{\threematrix}[1]{\left[\barr{ccc}#1\earr\right]}

	\newc{\mN}{\mathbf{N}}
	\newc{\mI}{\mathbf{I}}
	\newc{\mH}{\mathbf{H}}

	\newc{\mk}{\mathbf{k}}
	\newc{\mr}{\mathbf{r}}

	\newc{\DIAGM}[2]{\left[\barr{ccc}#1&0\ldots&0\\
		\vdots&\ddots&\vdots\\0&\ldots0&#2\earr \right]}
	\newc{\DiagM}[2]{\mbox{diag}\left[#1
		\cdots #2 \right]}
	\newc{\vVEC}[2]{\left[\barr{c}#1\\
		\vdots\\#2\earr \right]}
	\newc{\hVEC}[2]{\left[#1
		\cdots #2 \right]}
	
	\newc{\mq}{\mathbf{q}}

	\newc{\msys}[1]{\left\{\barr{l}#1\earr
		\right.}
	\newc{\msysa}[1]{\left\{\barr{ll}#1\earr
		\right.}
	
	\newc{\bbM}{\mathbb{M}}
	\newc{\mat}[1]{\left[\barr{cc}#1\earr\right]}
	
	\newc{\me}{\mathbf{e}}

	\newc{\vecc}[2]{\left[\barr{cc}#1\\#2\earr\right]}
	\newc{\mL}{\mathbb{L}}
	
	\newc{\cO}{{\cal O}}
	
	\newc{\cM}{{\cal M}}
	
	\newc{\myega }{\eg_0(R)}
	\newc{\myeg}{\eg_1(\eg_*)}
	\newc{\myegp}{\hat{\eg}_1(\eg_*)}

\newc{\diagA}{\mathbb{A}_d}
\newc{\mBB}{\mathbb{B}}
\newc{\MLT}[1]{{\cal M}_{lt}(\Og,#1)}
\newc{\ALT}[1]{{\cal A}_{l}(\Og,#1)}
\newc{\mM}{\mathbb{M}}

\newc{\diag}[1]{\mbox{diag}(#1)}
\newc{\off}[1]{\mbox{offdiag}(#1)}
\newc{\mT}{\mathbb{T}}

\usepackage{amssymb}

\newc{\idmu}[2]{{\dspl\int_{#1}~#2~d\mu}}
\newc{\idllg}[2]{{\dspl\int_{#1}~#2~d\llg}}

\begin{document}

	\vspace*{-.8in}
	\begin{center} {\LARGE\em On the smallness of  mean oscillations on metric-measure spaces and applications}
		
	\end{center}

	\vspace{.1in}
	
	\begin{center}

		{\sc Dung Le}{\footnote {Department of Mathematics, University of
				Texas at San
				Antonio, One UTSA Circle, San Antonio, TX 78249. {\tt Email: Dung.Le@utsa.edu}
				}}

	\end{center}

	\begin{abstract}
		It will be established that the mean oscillation of a function on a metric-measure space $X\times Y$ will be small if its mean oscillation on $X$ is small and some simple information on its (partial $Y$) upper-gradient is given. Applications to the regularity and global existence of bounded solutions to strongly coupled elliptic/parabolic systems on thin domains are also considered. \end{abstract}
	
\section{Introduction} In this paper we investigate the following property

\bdes\item[BMOsmall)] The BMO norm $\|u\|_{BMO(B_R)}$ is small if $R$ is sufficiently small.
\edes
Here, $B_R$ denotes a ball in $\RR^N$, or a metric space, with radius $R$. In \cite{dlebook, dlebook1}, we elucidate the pivotal role played by the smallness of BMO (bounded mean oscillation) norms in small balls within the examination of the regularity (and, in some cases, the existence) of weak solutions to strongly coupled elliptic/parabolic partial differential systems \cite{Am2, Gius} on $\RR^N$. This paper scrutinizes BMOsmall) in detail and introduces novel theoretical functional results that are applicable to solutions of partial differential systems on thin domains.

Moreover, we present broader theoretical findings in metric-measure spaces, building upon the theories presented in \cite{Haj, SmeetP}. Under general settings and employing straightforward arguments, we establish that the mean oscillation of a function on $X \times Y$ will be small if its mean oscillation on $X$ is small, provided there is some basic information on its (partial $Y$) upper-gradient. Several consequences ensue when these assumptions are easily verifiable, especially when the measure and metric on $X$ are typically related. These results are detailed in \refsec{smallBMOmetric}.

We now consider $X, Y$ as Euclidean spaces $\RR^N, \RR^K$ with their usual Lebesgue measures. In \refsec{euclid}, we apply the abstract results from \refsec{smallBMOmetric} to derive a condition on spatial partial derivatives for verifying BMOsmall). This condition is subsequently applied to bounded weak solutions of regular elliptic/parabolic second-order partial differential systems. In conjunction with the findings in \cite{GiaS, Gius}, we demonstrate that bounded weak solutions exhibit Hölder continuity if the domain is sufficiently thin. Assuming further that $\ma$ satisfies a spectral-gap condition, we establish that bounded weak solutions are classical and {\em exist globally}, as expounded in the theories of \cite{dlebook, dlebook1}.

It is worth noting that the exploration of scalar reaction-diffusion equations on 'thin domains' within the context of $\RR^2$ has been a subject of extensive research in the existing literature. This investigation aims to unravel the dynamics inherent in the (global) solutions of such equations. The global existence result on much more general thin domains in $\RR^N$, $N>2$, may serve as the inception of this research trajectory by demonstrating that strongly coupled parabolic systems define global dynamical systems of $m$ equations ($m>1$) within spaces of smooth functions.

	\section{The smallness of  mean oscillations on metric-measure spaces }\label{smallBMOmetric}\eqnoset
In this section, we study the general BMOsmall) property in the metric-measure spaces.

\subsection{A theoretic functional result} \label{theoretic}
We follow the notations in \cite{Haj}. Let $X$ be a metric-measure space with measure $\mu$. We assume that every pairs of points $x,y$ in $X$ can be connected by a rectifiable curve $\cg:[a,b]\to X$. A Borel (nonnegative) function $g$ is called an upper-gradient of $u:X\to\RR$ if $$|u(\cg(a))-u(\cg(b))|\le \int_\cg g \mbox{ for every rectifiable curve $\cg:[a,b]\to X$. }$$

Let $Y$ be another metric-measure space with measure $\llg$. Consider a function  $u:X\times Y\to \RR$ which is measurable with respect to the measure $\mu$. Define $U:Y\to\RR$ by $$U(y)=\idmu{X}{u(x,y)}, \quad y\in Y.$$ 

Suppose that $g$ is an ({\em partial}) upper-gradient of $u(x,\cdot)$ in $Y$ then it is easy to see that an upper-gradient of $U$ in $Y$ is $$G=\idmu{X}{g}.$$ Indeed, by Fubini's theorem, if $\cg:[a,b]\to Y$ is a rectifiable curve in $Y$ then
$$|U(\cg(a))-U(\cg(b))|\le \idmu{X}{|u(\cg(a))-u(\cg(b))|}\le \idmu{X}{\int_\cg g}=\int_\cg\idmu{X}{g}=\int_\cg G.$$

Following \cite[section 8]{Haj}, for some fixed  $b,s>0$, $\sg\ge1$ and fixed ball $B_{R_0}$ in $X$ we say that a measure  
\beqno{muR} \mbox{$\mu$ satisfies the condition $V(B_{\sg R_0},s,b)$ if } \mu(B_r(x))\ge br^s,\quad \forall B_r(x)\subset B_{\sg R_0} .\eeq  We then have the following Poincar\'e inequality for $p_*,p\ge 1$ and $r>0$ (see \cite[Remark 8.8]{Haj}) 
\beqno{Pineqmetric} \left(\frac{1}{\mu(B_r)}\int_{B_r}|u-u_{B_r}|^{p_*}d\mu\right)^\frac{1}{p_*}\le \left(\frac{\mu(B_{\sg r})}{r^s}\right)^\frac{1}{p}r\left(\frac{1}{\llg(B_{\sg r})}\int_{B_{\sg r}}g^{p}d\mu\right)^\frac{1}{p}.    \eeq

Let $X,Y$ be metric-measure spaces with measures $\mu,\llg$ respectively. We assume that every pairs of points $x,y$ in $Y$ can be connected by a rectifiable curve $\cg:[a,b]\to Y$ and $\llg$ satisfies $V(B'_{\sg R_0},l,b)$ \mref{muR}. {\em Notice that we do NOT assume this property for $X,\mu$}.

Let $B_R$ (respectively $ B'_r$) be a cube  in $X$ (respectively $Y$).  In the sequel,  we estimate  the following quantity for an integrable  function $W$ on $B_R\times B'_r$ ($W_{B_R}=W_{B_R}(y_*)$ is the $\mu$-average of $W$ over $B_R\times\{y_*\}$,  $y_*$ is dropped for simplicity of notations in the sequel)
\beqno{m0}\frac{1}{\mu(B_R) \llg(B'_r)}\dspl{\int_{B'_r}}\idmu{B_R}{|W-W_{B_R}|}d\llg =   \frac{1}{\llg(B'_r)}\dspl{\int_{B'_r}}\frac{1}{\mu(B_R)} \idmu{B_R}{|W-W_{B_R}|}d\llg.\eeq

We have our first result in the following theorem. {\em We should remark that there is no hypothesis on the relation between the metrics and measures of $X,Y$.}

\btheo{KEYBMOthmmetric} Let $X,Y$ be separable metric-measure spaces with $\sg$-finite Borel measures $\mu,\llg$ respectively. We assume that 

\bdes\item[1)] every pairs of points $x,y$ in $Y$ can be connected by a rectifiable curve $\cg:[a,b]\to Y$,
\item [2)] $\llg$ is doubling and satisfies the condition $V(B'_{\sg R_0},l,b)$ \mref{muR} for some $b,l,R_0>0, \sg\ge1$. \edes

For some  $\eg>0$, let $B_R, B'_r$ be  balls  in $X,Y$ respectively and let $W$ be a integrable function on $B_R\times B'_r$. We also assume that  $\mu(B_R\times \{y_*\})=\mu(B_R)$ for all $y_*\in Y$, and
\bdes
\item [i)]   for some ({\em partial}) upper gradient $g$ in $Y$ of $W$ and some finite $p\ge1$ we have   
\beqno{DxNW1a}\int_{B'_{\sg r}\times B_R}g^p d\mu d\llg=\int_{B'_{\sg r}}\idmu{B_R}{g^p}d\llg\le \eg r^{l-p}\mu(B_R),\eeq

\item [ii)]   for {\bf some} $s_*\in Y$  such that  the center $s_*$ of $B'_r$ is a Lebesgue
point of the function \beqno{fdef}f(y_*)=\frac{1}{\mu(B_R)}\idmu{B_R}{|W(x,y_*)-W_{B_R}|},\quad y_*\in Y.\eeq
Suppose that $R,r$ are small such that
$|f(s_*)|\le \eg$.
\edes
Then, we can assert that for the cube $\mB_{R,r}=B_R\times B'_r(s_*)$ in $X\times Y$ (with complete product measure $\mu\times\llg$) $$\frac{1}{\mu\times\llg(\mB_{R,r})}\idmu{\mB_{R,r}}{|W-W_{\mB_{R,r}}|}d\llg\le C(\eg^\frac{1}{p}+\eg)\quad \mbox{for some constant $C$} .$$

\etheo

\bproof For  $f(y_*)$ as in \mref{fdef}, $y_*\in Y$, and $g(x,\cdot)$ be an upper-gradient of $W$ in $Y$ (for each fixed $x\in X$) it is easy to see that (by Fubini-Tonelli's theorem) an upper-gradient of $f$ in $Y$ is \beqno{Gg}G(y_*)=\frac{1}{\mu(B_R)}\idmu{B_R}{2g}.\eeq

Applying the Poincar\'e inequality \mref{Pineqmetric} for $u=f, s=l$ on $Y$ for $r=R$, $p_*=1$ and $p\ge 1$, we have 
\beqno{m1}\frac{1}{\llg(B'_R)}\int_{B'_R}|f-f_{B'_R}|d\llg\le \left(\frac{\llg(B'_{\sg R})}{R^l}\right)^\frac{1}{p}R\left(\frac{1}{\llg(B'_{\sg R})}\int_{B'_{\sg R}}G^{p}d\llg\right)^\frac{1}{p}.    \eeq

Therefore, by \mref{Gg}, the  quantity \mref{m0} is majorized by
\beqno{m1a} \left(\frac{\llg(B'_{\sg R})}{R^l}\right)^\frac{1}{p}R\left(\frac{1}{\llg(B'_{\sg R})}\int_{B'_{\sg R}}\left(\frac{1}{\mu(B_R)}\idmu{B_R}{2g}\right)^{p}d\llg\right)^\frac{1}{p}+|f_{B'_R}|.\eeq
By H\"older's inequality, $\frac{1}{p}+\frac{1}{q}=1$,
$$\left(\frac{1}{\mu(B_R)}\idmu{B_R}{2g}\right)^{p} \le C    \mu(B_R)^{\frac{p}{q}-p}\idmu{B_R}{g^p}.$$
So that, after some simple algebraic simplifications, the  quantity \mref{m0} is also majorized by
\beqno{keybound} CR^{1-\frac{l}{p}}\mu(B_R)^{-\frac{1}{p}}\left(\int_{B'_{\sg R}}\idmu{B_R}{g^p}d\llg\right)^\frac{1}{p}+|f_{B'_R}|.\eeq

As $\llg$ is doubling (so that the Lebesgue differentiation theorem holds, see \cite[Section 2.9]{Federer}) and the center $s_*$ of $B'_R$ is a Lebesgue
point for $f$. Then if $\eg>0$ and $|f(s_*)|\le \eg$ then for $R$ is sufficiently small   we have $|f_{B'_R}|\le \eg$.
Hence, if for such $R$ we have \mref{DxNW1a} 
then this and \mref{keybound} show that the  quantity \mref{m0} is also majorized by $\eg^\frac{1}{p}+\eg$.

Of course, writing $W_{\mB_{R,r}}=\frac{1}{\llg(B'_r)}\dspl{\int}_{B'_r}W_{B_R}d\llg$ and $\frac{1}{\mu\times\llg(\mB_{R,r})}\idmu{\mB_{R,r}}{|W-W_{\mB_R}|}d\llg$ as
$$\frac{1}{\mu(B_R)\llg(B'_r)}\left[ \idmu{\mB_{R,r}\cap\{W>W_{\mB_{R,r}}\}}{(W-W_{\mB_{R,r}})}d\llg-\idmu{\mB_{R,r}\cap\{W\le W_{\mB_{R,r}}\}}{(W-W_{\mB_{R,r}})}d\llg\right],$$
we obtain the assertion of our theorem.
\eproof

We note that the assumption that $\llg$ is doubling in 2) was only for the Lebesgue differentiation theorem holds only (see \cite[Section 2.9]{Federer} for other cases). The Poincar\'e inequality \mref{Pineqmetric} does not need this assumption. In fact, the assumption that $\llg$ is doubling also implies that $\llg$ satisfies the condition $V(B'_{\sg R_0},l,b)$ for some $b,l,R_0>0, \sg\ge1$(\cite[Lemma 14.6]{SmeetP}).

Until now, we have not assume 1)-2) for $X,\mu$. Assume 1)-2) here for $\mu$. Under more assumptions on the relation of $\mu,\llg$, we have the following consequence.
\bcoro{KEYBMOcorometric1} Let $X,Y,\mu,\llg$ be as in \reftheo{KEYBMOthmmetric}. Assume that $\mu$ satisfies the condition $V(B_{\sg R_0},s,b)$ \mref{muR} for some $b,s,R_0>0, \sg\ge1$. In particular, we assume \mref{muR}, $l\le s$, and for some $L\ge1$ we have $\llg(B'_{R^L})\le \mu(B_R)$ for all $0<R\le R_0$. Suppose that $R_0$ is small.

If we suppose that $\{|g(\cdot,y)|^p\}_{y\in B'_{\sg R_0}}$, with  $p=s\ge 1$, is {\em uniformly} $\mu$-integrable, then  for the cube $\mB_{R,r}=B_R\times B'_r$ in $X\times Y$ with $r=R^L$ and complete product measure $\mu\times\llg$ $$\frac{1}{\mu\times\llg(\mB_{R,r})}\idmu{\mB_{R,r}}{|W-W_{\mB_{R,r}}|}d\llg\le C(\eg^\frac{1}{p}+\eg)\quad \mbox{for some constant $C$} .$$

\ecoro

\bproof
For some constants $c_1,c_2>0$, by Jensen's inequality and \cite[ii) of Theorem 8.7]{Haj} for $\mu$, we have $$\exp\left(\frac{1}{\mu(B_R)}\idmu{B_R}{\left[c_1 b^\frac{1}{s}\frac{|W-W_{B_R}|}{\|g\|_{L^p(B_{\sg R})}}\right]} \right)\le \frac{1}{\mu(B_R)}\idmu{B_R}{\exp\left(c_1 b^\frac{1}{s}\frac{|W-W_{B_R}|}{\|g\|_{L^p(B_{\sg R})}}\right)} \le c_2,  $$ where $\|g\|_{L^p(B_{\sg R})}=\|g\|_{L^p(B_{\sg R},\mu)}$. This implies
$$\frac{1}{\mu(B_R)}\idmu{B_R}{|W-W_{B_R}|} \le \frac{\|g\|_{L^p(B_{\sg R})}}{c_1 b^\frac{1}{s}}\ln(c_2).$$

Since $\mu(B_{\sg R)}$ is small if $R$ is small, we have $\|g\|_{L^p(B_{\sg R})}\le \eg$. Hence, the above yields $|f(s_*)|\le\eg$. This is i) of \reftheo{KEYBMOthmmetric}.  

Also, as we suppose that $\{|g(\cdot,y)|^p\}_{y\in B'_{\sg R_0}}$, with  $p=s\ge 1$, is {\em uniformly} $\mu$-integrable, if $R$ is sufficiently small we take  $r=R^L\le R$  such that $\llg(B'_{\sg R^L})\le c\llg(B'_{R^L}) \le c\mu(B_R)$ then
$$\int_{B'_{\sg r}}\idmu{B_R}{g^{p}}d\llg\le \eg^p \llg(B'_{\sg r})\le \eg \llg(B'_{\sg R^L})\le \eg c\mu(B_R).$$
So, as $l-p\le0$ and $r$ is small, we can have $c\le r^{l-p}$ and then
$$\int_{B'_{\sg r}}\idmu{B_R}{g^{p}}d\llg\le \eg r^{l-p}\mu(B_R).$$
This is \mref{DxNW1a} of ii). \reftheo{KEYBMOthmmetric} then completes the proof. \eproof

\bcoro{KEYBMOcorometric2} Let $X,Y,\mu,\llg$ be as in \refcoro{KEYBMOcorometric1}. Assume that $\llg(B'_r)\sim r^l$. Then, we can drop the assumption $l\le s$ and omit $L$.

If we suppose that $\{|g(\cdot,y)|^p\}_{y\in B'_{\sg R_0}}$, with  $p=s\ge 1$, is {\em uniformly} $\mu$-integrable, then  for the cube $\mB_{R,r}=B_R\times B'_r$ in $X\times Y$ with complete product measure $\mu\times\llg$ $$\frac{1}{\mu\times\llg(\mB_{R,r})}\idmu{\mB_{R,r}}{|W-W_{\mB_{R,r}}|}d\llg\le C(\eg^\frac{1}{p}+\eg)\quad \mbox{for some constant $C$} .$$

\ecoro

\bproof
In fact, if $\llg(B'_r)\sim r^l$ then (as $p=s\ge1$)
$$\int_{B'_{\sg r}}\idmu{B_R}{g^{p}}d\llg\le \eg^p\llg(B'_{\sg r})\le c\eg r^l\le \eg r^{l-p}r^s\le \eg r^{l-p}\mu(B_R).$$ The proof is complete. \eproof

The abstract assumption that $\{|g(\cdot,y)|^p\}_{y\in B'_{\sg R_0}}$, with  $p=s\ge 1$, is {\em uniformly} $\mu$-integrable can  be verified in some applications. However, we  also have a  much simpler following 

\bcoro{KEYBMOcorometric3} The assertions of \refcoro{KEYBMOcorometric1} and \refcoro{KEYBMOcorometric2} holds if
$(x_*,s_*)$, the center of $B_R\times B'_{\sg r}$  is a Lebesgue point of $g$ in $X\times Y$ with measure $\mu\times\llg$. 
\ecoro

\bproof
Indeed, by the assumption, for any given $\eg>0$ if $r,R$ (and therefore $\mu(B_{R})$) are small in terms of $g(x_*,s_*)^p,\eg$ then
$$\int_{B'_{\sg r}}\idmu{B_R}{g^{p}}d\llg\le 2g(x_*,s_*)^p\mu(B_R)\llg(B'_{\sg r})\Rightarrow \int_{B'_{\sg r}}\idmu{B_R}{g^{p}}d\llg\le \eg^p\llg(B'_{\sg r})$$
and the assertion follows from the proof of \refcoro{KEYBMOcorometric2}. \eproof

For $p<p_*$, by H\"older's inequality and continuity of integral we have for small $r,R$
$$\int_{B'_{\sg r}}\idmu{B_R}{g^{p}}d\llg\le  \left(\int_{B'_{\sg r}}\idmu{B_R}{g^{p_*}}d\llg \right)^\frac{p}{p_*}\left(\llg(B'_{\sg r})\mu(B_R)\right)^{1-\frac{p}{p_*}} \le \eg \left(\llg(B'_{\sg r})\mu(B_R)\right)^{1-\frac{p}{p_*}}.$$

$$\int_{B'_{\sg r}}\idmu{B_R}{g^{p}}d\llg\le   \eg \frac{\left(\llg(B'_{sg r})\right)^{1-\frac{p}{p_*}}}{\mu(B_R)^\frac{p}{p_*}}\mu(B_R).$$

We want $\frac{\left(\llg(B'_{sg r})\right)^{1-\frac{p}{p_*}}}{\mu(B_R)^\frac{p}{p_*}} \le r^{l-p}$. 

If $\llg(B'_{ R})\mu(B_R)\le cR^{s+l}$ and $p_*<s+l$, then $(s+l)(1-\frac{p}{p_*})>s+l-p$. Because $R<1$

$$\int_{B'_{\sg r}}\idmu{B_R}{g^{p}}d\llg\le   \eg R^{s+l-p}\le \eg R^{l-p}\mu(B_R)$$

{\em If we want to apply this result to functions on Euclidean spaces then we must have the continuity $\int_{B'_{\sg r}}\idmu{B_R}{g^{p_*}}d\llg\le \eg$ when $r,R$ small. This not true for solutions to elliptic systems although we have higher integrability of gradients. However, we can argue as in \cite[b) of Theorem 5.4]{keylist}.  }

$(s+l)(1-\frac{p}{p_*})-(s+l-p)=(s+l)(\frac{p}{s+l}-\frac{p}{p_*})$

Consider again \mref{m1a}. Simplify the first term and write it as $$\frac{R^{1-\frac{l}{p}}}{\mu(B_R)}\left(\int_{B'_{\sg R}}\left(\idmu{B_R}{2g}\right)^{p}d\llg\right)^\frac{1}{p}+|f_{B'_R}|.$$
Let $M$ be a nonnegative measurable function on $X\times Y$. Using H\"older's inequality, $\frac{1}{p}+\frac{1}{q}=1$, we have
$$\left(\idmu{B_R}{2g}\right)^{p} = \left(\idmu{B_R}{Mg2M^{-1}}\right)^{p}\le C    \|M^{-1}\|_{L^q(B_R,\mu)}^{p}\idmu{B_R}{M^pg^p}.$$
So that,
$$\left(\int_{B'_{\sg R}}\left(\idmu{B_R}{2g}\right)^{p}d\llg\right)^\frac{1}{p}\le $$

\subsection{A variance of \reftheo{KEYBMOthmmetric}} \label{var}
Let us go back to the proof of \reftheo{KEYBMOthmmetric} and begin with \mref{m0}. 

As before, we consider
$f(y_*)=\frac{1}{\mu(B_R)}\idmu{B_R}{|W(x,y_*)-W_{B_R}|}$, $y_*\in Y$, and let $g(x,\cdot)$ be an upper-gradient of $W$ in $Y$ (for each fixed $x\in X$) such that $g$ is integrable on $X\times Y$ with complete product measure $\mu\times\llg$. Importantly, we assume that $\mu(B_R\times \{y_*\})=\mu(B_R)$ for all $y_*\in Y$. An upper-gradient of $f$ in $Y$ is $$G(y_*)=\frac{1}{\mu(B_R)}\idmu{B_R}{2g}.$$

Instead of applying the Poincar\'e inequality \mref{Pineqmetric}, we have 
$$|f(y_*)-f(s_*)|\le \int_\cg G=\int_\cg\frac{1}{\mu(B_R)}\idmu{B_R}{2g}.$$

Suppose that $y_*,s_*$ are connected by a curve $\cg:[0,1]\to Y$, and $|\cg'|\le C d_Y(y_*,s_*)^\sg$, the metric in $Y$ is $d_Y$ and $\sg\ge 1$. Then we have ($g$ is an upper-gradient w.r.t this family of $\cg$)
$$|f(y_*)|\le \frac{d_Y(y_*,s_*)^\sg}{\mu(B_R)}\idmu{B_R}{2g}+|f(s_*)|.$$

Therefore, if $M$ is a positive (a.e) function on $Y$ then $I:=\frac{1}{\llg(B'_r)}\dspl{\int_{B'_R}}\frac{1}{\mu(B_R)} \idmu{B_R}{|W-W_{B_R}|}d\llg$ is estimated by ($I$ is the quantity in \mref{m0} and $s_*$ is the center of $B'_r$)
$$  \frac{2}{\llg(B'_r)\mu(B_R)}\dspl{\int_{B'_r}}M^{-1}(y_*) d_Y(y_*,s_*)^\sg\idmu{B_R}{M(y_*)g} d\llg+|f(s_*)|. $$

For $\frac{1}{p}+\frac{1}{q}=1$, by H\"older's inequality on $(Y,\llg)$
$$\dspl{\int_{B'_r}}M^{-1}(y_*) d_Y(y_*,s_*)^\sg\idmu{B_R}{M(y_*)g}\le \|M^{-1}(y_*) d_Y(y_*,s_*)\|_{L^q(B'_r,\llg)}\left(\dspl{\int_{B'_r}}\left[\idmu{B_R}{M(y_*)g}\right]^p d\llg\right)^\frac{1}{p} $$

Of course, $\left[\idmu{B_R}{M(y_*)g}\right]^p\le \mu(B_R)^\frac{p}{q}\idmu{B_R}{M^pg^p}$, so that we  derive 
\beqno{I0}I\le \frac{2}{\llg(B'_r)\mu(B_R)^{\frac{1}{p}}}\|M^{-1}(y_*) d_Y(y_*,s_*)^\sg\|_{L^q(B'_R,\llg)}\left(\dspl{\int_{B'_r}}\idmu{B_R}{M^p(y_*)g^p} d\llg\right)^\frac{1}{p}+|f(s_*)|. \eeq

If we take $M(y_*)=d_Y(y_*,s_*)^{\sg-\frac{s}{p}}$ {\em for some $p$ such that $\sg p>s$}. Then, $|M^{-1}d_Y(y_*,s_*)^\sg|^q=d_Y(y_*,s_*)^\frac{sq}{p}$ and hence $\|M^{-1}(y_*) d_Y(y_*,s_*)^\sg\|_{L^q(B'_R,\llg)}=\left( \int_{B'_r} d_Y(y_*,s_*)^\frac{sq}{p}d\llg\right)^\frac{1}{q}\le c r^{\frac{s}{p}+\frac{L}{q}}$ if we assume $\llg(B'_r)\le br^L$ for some $L\le l$.

The first term on the right hand side of \mref{I0} is bounded by (if $\llg(B'_r)\ge b r^l$)
$$c\frac{r^{\frac{s}{p}+\frac{L}{q}}}{\llg(B'_r)\mu(B_R)^{\frac{1}{p}}}\left(\dspl{\int_{B'_r}}d_Y(y_*,s_*)^{\sg p-s}\idmu{B_R}{g^p} d\llg\right)^\frac{1}{p}\le C\frac{r^{\frac{s}{p}+\frac{L}{q}-l}}{\mu(B_R)^{\frac{1}{p}}}\left(\dspl{\int_{B'_r}}d_Y(y_*,s_*)^{\sg p-s}\idmu{B_R}{g^p} d\llg\right)^\frac{1}{p}
$$
which will be less than $\eg^\frac{1}{p}$ for any given $\eg>0$ if

$$\dspl{\int_{B'_r}}d_Y(y_*,s_*)^{\sg p-s}\idmu{B_R}{g^p} d\llg\le \eg r^{-p(\frac{s}{p}+\frac{L}{q}-l)}\mu(B_R)=\eg r^{-p(\frac{L}{q}-l)-s}\mu(B_R).$$

Notice that we do not assume that $s_*$ is a Lebesgue point in $Y,\llg$ of any function. In addition, we need not the Poincar\'e's inequality \mref{Pineqmetric} (although it is available) as in the proof of \reftheo{KEYBMOthmmetric}. We gather our assumptions and assertions in  the following 
\btheo{KEYBMOthmmetric1} Let $X,Y$ be separable metric-measure spaces with $\sg$-finite Borel measures $\mu,\llg$ respectively. We assume that 

\bdes\item[1)] every pairs of points $s_*,y_*$ in $Y$ can be connected by a rectifiable curve $\cg:[0,1]\to Y$, and $|\cg'|\le C d_Y(y_*,s_*)^\sg$ for some constants $C>0$ and $\sg\ge1$,
\item [2)] $\llg$ satisfies the condition $b'r^L\ge \llg(B_r)\ge br^l$ for all $B_r\subset B_{R_0} $  for some constants $b,b'>0$, $L\le l,$ and $0<r<R_0$ small (see \mref{muR}). \edes

For some  $\eg>0$, let $B_R, B'_r(s_*)$ be  balls  in $X,Y$ respectively, and $W$ be an integrable  function on $B_R\times B'_r(s_*)$. We also assume that  $\mu(B_R\times \{y_*\})=\mu(B_R)$ for all $y_*\in Y$, and
\bdes
\item [i)]   for some ({\em partial}) upper gradient $g$ in $Y$ of $W$, $g$ is an upper-gradient w.r.t this family of $\cg$ in 1), and {\em some finite $\sg p>s$} we have  
\beqno{DxNW1b}\dspl{\int_{B'_r(s_*)}}d_Y(y_*,s_*)^{\sg p-s}\idmu{B_R}{g^p} d\llg\le \eg r^{-p(\frac{L}{q}-l)-s}\mu(B_R), \quad\mbox{ ($\frac{1}{p} +\frac{1}{q}=1$)}\eeq

\item [ii)]   for {\bf some} $s_*\in Y$ and $R,r$ are small such that  $$\frac{1}{\mu(B_R)}\idmu{B_R}{|W(x,s_*)-W_{B_R}|}\le \eg.$$
\edes
Then, we can assert that for the cube $\mB_{R,r}=B_R\times B'_r(s_*)$ in $X\times Y$ (with complete product measure $\mu\times\llg$) $$\frac{1}{\mu\times\llg(\mB_{R,r})}\idmu{\mB_{R,r}}{|W-W_{\mB_{R,r}}|}d\llg\le C(\eg^\frac{1}{p}+\eg)\quad \mbox{for some constant $C$} .$$

\etheo

\brem{pintegrabily} The assumption that $d_Y(y_*,s_*)^{\sg p-s}g^p$ is (locally) integrable on $X\times Y$ is also very important so that the integral in \mref{DxNW1b} makes sense. \erem

\brem{BMOthmmetricrem} Since $d_Y(y_*,s_*)\le r$, the condition \mref{DxNW1b} is  weaker than \mref{DxNW1a} in ii) of \reftheo{KEYBMOthmmetric} when $\sg=1$ and $L= l$ (so that $-p(\frac{L}{q}-l)= l$). Indeed, as $d_Y(y_*,s_*)\le r<1$ and $p>s$ 
$$\dspl{\int_{B'_r}}d_Y(y_*,s_*)^{p-s}\idmu{B_R}{g^p} d\llg\le \overbrace{r^{p-s}\int_{B'_{\sg r}}\idmu{B_R}{g^p}d\llg\le \eg r^{l-s}\mu(B_R)}^{\mbox{\mref{DxNW1a}}}= \eg r^{-p(\frac{L}{q}-l)-s}\mu(B_R).$$

When $L\le l$, we cannot compare the two conditions. We have $-p(\frac{L}{q}-l)\ge l$ so that 
$$\dspl{\int_{B'_r}}d_Y(y_*,s_*)^{p-s}\idmu{B_R}{g^p} d\llg\le r^{p-s}\int_{B'_{\sg r}}\idmu{B_R}{g^p}d\llg\le  \eg r^{-p(\frac{L}{q}-l)-s}\mu(B_R)\le \eg r^{l-s}\mu(B_R).$$
\erem

\brem{psrem}If $\mu(B_R)\sim R^s$, $L=l$, and $r=R$ then \mref{DxNW1b} is simply 
$$\dspl{\int_{B'_R(s_*)}}d_Y(y_*,s_*)^{\sg p-s}\idmu{B_R}{g^p} d\llg\le \eg R^{l},\quad \mbox{for some $p$ such that }\sg p>s.$$
\erem

\section{Euclidean spaces and applications}\label{euclid}\eqnoset

We now taking $X,Y$ to be Euclidean spaces $\RR^N,\RR^K$ with their usual Lebesgue measures. All the assumptions 1)-2) of \reftheo{KEYBMOthmmetric} are verified but we choose to state the assumption i) a bit stronger according to \mref{DxNW1a} of \reftheo{KEYBMOthmmetric} (see \refrem{BMOthmmetricrem}) in the following

\btheo{KEYBMOthm} Let $K\ge1$ and $\Og$ be a domain in $\RR^{N+K}$ and  $\mB_R$ be a cube in $\Og$ with sides parallel to the axes of $\RR^{N+K}$. We write $\RR^{N+K}=\{(x,y_*)\,:\, x\in\RR^N, y_*\in\RR^K\}$.   If for small $R_0,\eg>0$

\bdes
\item [i)] $D_{y_*}W \in L^{p,N+K-p}(\Og)$  for some finite $p\ge1$ and for all $0<R\le R_0$
\beqno{DxNW1}\dspl{\dspl{\int_{B'_R(s_*)}}} \iidx{B_R}{| D_{x_{*}}W(x,y_*)|^p}dy_*\le \eg R^{N+K-p}, \quad B_R\subset\RR^{N},\;B'_R(s_*)\subset\RR^{K};\eeq

\item [ii)]   for {\bf some} $s_*\in\RR^K$  $$f(y_*)=\frac{1}{R^N}\iidx{B_R}{|W(x,y_*)-W_{B_R}|},\quad y_*\in Y.$$
Suppose that $R$ is small such that
$|f(s_*)|\le \eg$ and $|f_{B'_{R}}|\le \eg$.
\edes
Then, for $R_0,\eg$ as in \mref{DxNW1}, $\mB_R=B'_R\times B_R(s_*)$, and for some constant $C$
\beqno{notBMO}\frac{1}{R^{N+K}}\int_{\mB_R}{|W-W_{\mB_{R}}|}dxdy_*\le C(\eg^\frac{1}{p}+\eg) .\eeq

\etheo

It is important to note that \mref{notBMO} does NOT say that the BMO norm of a function $W$ on a cube $\mB_R$  with sides parallel to the axes of $\RR^{N+K}$ will be small in $R$ if its BMO norms in $N$ dimensional cubes are also small and \mref{DxNW1} holds for some small $\eg>0$ and $p\ge1$ (a simple information on the $\RR^K$-direction derivatives). 

\brem{grem}
Of course, the upper-gradient $g$ in \reftheo{KEYBMOthmmetric} is exactly the partial derivative $D_{x_*}W$. If we use \reftheo{KEYBMOthmmetric1} and  $g$ is an upper-gradient w.r.t the family of rectifiable curve $\cg$ satisfying 1) of \reftheo{KEYBMOthmmetric1}, then $g$ is $D_{|x_*|^{\sg-1}x_*}W$, $\sg\ge1$.
\erem

However, in applications to {\em bounded} weak solutions to strongly coupled partial differential systems (see \cite{Gius,GiaS}), the assertion \mref{notBMO} is enough, if $\eg$ is sufficiently small, to conclude that they are H\"older continuous and we also have $\|W\|_{BMO(\mB_{R_0})}\le C(\eg^\frac{1}{p}+\eg)$.

Consider regular elliptic/parabolic second order partial differential systems
\beqno{exsyszpara}-\Div(\ma(W)DW)=F(W) \mbox{ or }W_t=\Div(\ma(W)DW)+F(W)\eeq
on bounded domain $\Og$ of $\RR^N$ or $\Og\times(0,T)$ with homogeneous Dirichlet or Neumann boundary conditions. Here, $W$ is an unknown vector in $\RR^m$, $m\ge 2$. As usual, $\ma(W), F(W)$ are continuous $m\times m$ matrix or $m$ vector and  for some positive constants $\llg_0,\LLg_0$ we assume the ellipticity condition
\beqno{normalellcond}\LLg_0|DW|^2\ge \myprod{\ma(W)DW,DW}\ge \llg_0|DW|^2.\eeq

We first discuss the elliptic case. 
We can use the theorem  on thin domains ($N=p=2$ and thin in the $K$-direction). We can let $\Og=S\times B'_\dg$ for small $\dg$ and $S\subset\RR^2$, a bounded domain, and assume the Neumann boundary condition $D_{x_*}W=0$ on $S$. Indeed, 
$$\iidx{B_R}{|D_{x_*} W|^2}=0 \mbox{ for $B_R\subset S$} \Rightarrow\lim_{R\to0}\frac{1}{R^K}\int_{B'_R}\iidx{B_R}{|D_{x_*} W|^2} dx_*=0$$ so that when $R$ is small (or $\dg$ is small) so that  the following holds for any cube in $\Og$
$$\frac{1}{R^K}\int_{B'_R}\iidx{B_R}{|D_{x_*} W|^2} dx_*\le \eg \Rightarrow \int_{B'_R}\iidx{B_R}{|D_{x_*} W|^2} dx_*\le \eg R^K=\eg R^{K+N-2}. $$

This is i) of \reftheo{KEYBMOthm}. The assumption ii) is obvious by Poincar\'e's inequality on $\RR^2$ as we are considering bounded $W$ and restrict \mref{exsyszpara} to 2-dimensional sections of $\Og$.

The parabolic case is similar as the above argument involves only with the spatial partial derivative $D_{x_*}$ in the $\RR^K$ direction (we then obtain the BMOsmall) in space-time and it is sufficient to use \cite{GiaS}). The assumption ii) can be established by following the planar argument in \cite{dleJMAA}.

Under \mref{normalellcond}, we now obtain the H\"older continuity of $W$ as in \cite[Theorem 3.1]{GiaS} by an argument based on a comparison with solutions to freezing coefficient systems in \cite[Proposition 1.2]{GiaS}. 
We then assert the following

\bcoro{continuitycoro} Consider a domain in $\RR^N$, $N>2$, $\Og=S\times B'_\dg$ for small $\dg>0$ and $S\subset\RR^2$, a bounded domain. 

Then, bounded weak solutions of regular elliptic/parabolic second order partial differential systems
$$-\Div(\ma(W)DW)=F(W) \mbox{ or }W_t=\Div(\ma(W)DW)+F(W)$$ satisfying \mref{normalellcond} are H\"older continuous if they satisfy homogeneous Neumann boundary condition on $S$ and $\dg$ is sufficiently small (i.e. $\Og$ is sufficiently thin).
\ecoro

Note that we cannot use \refcoro{KEYBMOcorometric3} as $R$ must be small in terms of $g(x_*,s_*)=D_{x_*}W$ and $\eg$. Therefore, combining with the partial H\"older regularity results, we can only assert that bounded weak solutions to elliptic/parabolic cross-diffusion systems are H\"older continuous almost everywhere as in \cite{GiaS}. This is different with scalar equations ({\em see also \cite[p. 82]{SmeetP}} for possible generalizations to more abstract settings).

On the other hand, suppose that there are $\llg(W),\LLg(W)>0$ such that $$\llg(W)|\zeta|^2\le \myprod{\ma(W)\zeta,\zeta}\le \LLg(W)|\zeta|^2\quad\forall\zeta\in\RR^m\times\RR^N.$$ Let $\nu_*=\sup_{W\in\RR^m}\frac{\llg(W)}{\LLg(W)}$. We say that 
$\ma(W)$ verifies the {\em spectral-gap condition} if  $\nu_*>1-2/N$.  Note that is condition is void if $m=1$.

Once $W$ is H\"older continuous, the BMOsmall) is satisfied and we can invoke the theory in \cite{dlebook,dlebook1} (where the so called strong/weak Gagliardo-Nirenber inequalities with BMO norm are developed) to assert that $W$ exists globally if it satisfies BMOsmall). We then have
\bcoro{existencecoro} Consider a domain in $\RR^N$, $N>2$, $\Og=S\times B'_\dg$ for small $\dg>0$ and $S\subset\RR^2$, a bounded domain. Assume that 
$\ma(W)$ verifies \mref{normalellcond} and {\em the spectral-gap condition}. 

Then, bounded weak solutions of regular parabolic second order partial differential systems
$$W_t=\Div(\ma(W)DW)+F(W)$$  are classical and exist {\em globally} if they satisfy homogeneous Neumann boundary condition on $S$ and $\dg$ is sufficiently small (i.e. $\Og$ is sufficiently thin).

\ecoro

Finally,
we remark that the exploration of scalar reaction-diffusion equations on 'thin domains' in $\RR^2$ has been a focal point in existing literature, seeking to comprehend the intricate dynamics embedded in their (global) solutions.  \refcoro{existencecoro} serves as the inception of this research trajectory by demonstrating that strongly coupled parabolic systems define global dynamical systems ($m>1$) within spaces of smooth functions.

\bibliographystyle{plain}

\end{document}